\documentclass[letterpaper, 10 pt, conference]{ieeeconf}  

\IEEEoverridecommandlockouts                              
\usepackage{algorithmic}
\usepackage{epstopdf}
\usepackage{float}
\usepackage{amsmath,amsfonts}
\usepackage{graphicx}

\title{\LARGE \bf
Stochastic Finite Volume Method for Uncertainty Management in \\ Gas Pipeline Network Flows
}
\author{Saif R. Kazi$^\star$, Sidhant Misra$^\star$, Svetlana Tokareva$^\star$, Kaarthik Sundar$^\dag$, and Anatoly Zlotnik$^\star$
\thanks{This project was supported by the LDRD program project ``Stochastic Finite Volume Method for Robust Optimization of Nonlinear Flows'' at Los Alamos National Laboratory.  Research conducted at Los Alamos National Laboratory is done under the auspices of the National Nuclear Security Administration of the U.S. Department of Energy under Contract No. 89233218CNA000001. Report No. LA-UR-24-22647.}
\thanks{$^\star$\{skazi,sidhant,tokareva,azlotnik\}@lanl.gov, \,\, Applied Mathematics \& Plasma Physics, Los Alamos National Laboratory, Los Alamos, NM 87545}
\thanks{$^\dag$kaarthik@lanl.gov \,\, Information Systems \& Modeling, Los Alamos National Laboratory, Los Alamos, NM 87545}}

\begin{document}

\maketitle
\thispagestyle{empty}
\pagestyle{empty}

\begin{abstract}
Natural gas consumption by users of pipeline networks is subject to increasing uncertainty that originates from the intermittent nature of electric power loads serviced by gas-fired generators.  To enable computationally efficient optimization of gas network flows subject to uncertainty, we develop a finite volume representation of stochastic solutions of hyperbolic partial differential equation (PDE) systems on graph-connected domains with nodal coupling and boundary conditions. The representation is used to express the physical constraints in stochastic optimization problems for gas flow allocation subject to uncertain parameters. The method is based on the stochastic finite volume approach that was recently developed for uncertainty quantification in transient flows represented by hyperbolic PDEs on graphs. In this study, we develop optimization formulations for steady-state gas flow over actuated transport networks subject to probabilistic constraints.  In addition to the distributions for the physical solutions, we examine the dual variables that are produced by way of the optimization, and interpret them as price distributions that quantify the financial volatility that arises through demand uncertainty modeled in an optimization-driven gas market mechanism.  We demonstrate the computation and distributional analysis using a single-pipe example and a small test network.
\end{abstract}

\section{Introduction} \label{sec:intro}

The increasing use of wind and solar electricity production concurrently with growing reliance on gas-fired power generation for both base load and to compensate for intermittent renewables results in more volatile and uncertain demands on natural gas delivered by pipelines \cite{yang2017effect}.  This issue has been addressed in part by developing probabilistic scheduling methods based on power grid physics and operations for calibrated planning of gas-fired electricity generation and its fuel consumption \cite{roald2017chance}.  For instance, stochastic optimization methods were developed based on sampling demand and price scenarios associated with wind uncertainty \cite{zhao2016unit}.  These methods take advantage of the understanding that power grid physics are approximated well for operations scheduling using direct current (DC) power flow modeling.  In contrast, it has been well known for decades that probabilistically constrained optimization for gas pipeline network flow allocation is highly challenging because of nonlinearity \cite{levary1980natural}.

The notion of robust optimization involves the formulation of mathematical programs to determine optimal decisions that are feasible under a range of parameter values \cite{sahinidis2004opportunities,bertsimas2011theory}.  Formulations based on scenario sampling have been developed for managing demand uncertainty in operations of steady and transient gas network flows \cite{behrooz2016managing,zavala2014stochastic}, including recent approaches that account for composition uncertainty \cite{liu2020dynamic}.  Recent formulations have proposed integrated optimization formulations to jointly manage day-ahead uncertainty in a power grid and interconnected gas pipeline, including compensation for interval \cite{roald2020uncertainty} and intertemporal \cite{kazi2024intertemporal} variability.  Despite these advances, the more fundamental steady-state optimal gas flow under demand uncertainty remains challenging to analyze and scale to large networks.

The interconnections between physical and financial management of the power grid and gas pipelines in the presence of fuel price volatility and renewable generation uncertainty raise complex economic issues \cite{geng2021network,wang2022volatility}.  Whereas economics are connected to operations in the power grid through well-established optimization-driven wholesale electricity markets \cite{schweppe2013spot}, such mechanisms for gas pipeline management have been explored much more recently \cite{pepper2012implementation}.  Deterministic optimization-driven market mechanisms for steady-state and dynamic hourly scheduling and pricing of natural gas delivery have been developed  \cite{rudkevich17hicss,zlotnik2019optimal}, and such formulations were shown to scale well to realistic systems \cite{zlotnik2017economic}.  However, the extension of market mechanisms for gas pipelines to stochastic formulations have not been examined to date. It is well known that optimization of gas pipeline flows under uncertainty is highly complex, because of the substantial nonlinearity and computational complexity involved \cite{levary1980natural,behrooz2016managing}.  The challenge of modeling the effect of intertemporal uncertainty in time-dependent simulation of gas pipelines has to an extent been met by a new stochastic finite volume (SFV) method for uncertainty quantification (UQ) in gas pipeline flows \cite{tokareva2024stochastic}, which was shown to model the propagation of uncertainty in associated initial boundary value problems.  Moreover, finite volume methods have been used in financial analysis of markets that do not involve underlying physics, such as options pricing \cite{huang2006fitted}. 

In this paper, we adapt the SFV representation for solutions of the nonlinear gas flow equations with uncertain boundary conditions to optimize the steady-state flow allocation in a pipeline system subject to uncertain demands.  We propose this type of uncertainty modeling because of its computationally efficient representation of the probability densities of state variables throughout a pipeline system, and in particular its amenability to higher-order reconstructions in representing probabilistic inequality constraints.  The primary contribution of our study is a chance-constrained optimization formulation that extends the problem of steady-state flow allocation for a gas pipeline system \cite{rudkevich17hicss}, resulting in a single solution for compressor configuration and optimized flows that extremizes an operational and/or economic objective in expectation.   By utilizing the proven method of chance-constraints, our pipeline demand uncertainty management approach provides probabilistic robustness guarantees that can be calibrated based on risk tolerance.  Furthermore, when the objective function is interpreted as an expected economic value created for users of the pipeline system, we examine the dual variables of the optimization solution and use them to reconstruct functions that we interpret as the probability distributions of locational values of gas.

The paper is structured as follows.  We describe gas pipeline flow modeling and the steady-state optimal gas flow optimization formulation in Section \ref{sec:ssogf}.  This is extended to the chance-constrained stochastic steady-state optimal gas flow formulation that accounts for uncertain loads using the SFV method in Section \ref{sec:ccssogf}.  In Section \ref{sec:examples}, we describe our computational implementation and two small interpretable case studies, then briefly propose an economic interpretation of the results, and finally conclude in Section \ref{sec:conc}.

\section{Steady-State Optimal Gas Flow Problem}  \label{sec:ssogf} 

We formulate a general robust optimization problem for steady-state gas flow, in which the objective minimizes average gas compressor power (for maximal operating efficiency) and maximizes the expected economic value for users of the system in the form of revenue collected by the pipeline operator.  The formulation admits situations in which gas consumptions at demand nodes are specified, as well as problems where the gas delivery to one or more nodes is an optimized variable.  In the latter case, price-quantity bids for withdrawal (by consumers) and injection (by suppliers) are provided, and we suppose that a market authority solves the steady-state flow allocation problem subject to uncertain flows that are unknown but within a specified range.  

\subsection{Gas Flow Modeling}  \label{sec:gas_modeling}

We consider a minimal turbulent flow approximation relating pressure to mass flow for compressible gas in a pipe.  The flow of gas through a pipe in the turbulent regime can be described by a simplification of the Euler equations  \cite{herty2010new}, 
\begin{align}
  &\partial_t \rho + \partial_x \varphi = 0, \,\,\, \partial_t \varphi + \partial_x p  = -\frac{\lambda}{2D} \frac{\varphi|\varphi|}{\rho}, \label{eq:euler:1}
\end{align}
where $\varphi$ is per-area mass flux (kg/m$^2$/s) of the gas, $p$ is pressure (Pa), and $\rho$ is gas density (kg/m$^3$), which are defined on a domain $x\in[0,L]$ at each time $t\in[0,T]$.  The parameters that determine flow capacity are the pipe-specific non-dimensional friction factor $\lambda$, pipe diameter $D$, and pipe length $L$, and the wave propagation speed $a$ that is a property of the gas that depends on composition and temperature.  Recall that the per-area mass flux is $\varphi=\rho\!\cdot \!u$, and is related to the total flow $\phi$ through a pipe by $\varphi=\phi/A$ where $A=\pi D^2/4$ is the cross-sectional area of the pipe.  We assume that the wave speed $a$ is uniform system-wide, and that gas pressure $p$ and density $\rho$ are related by the ideal gas law $p = a^2 \rho$.  We simplify this way to focus on probabilistic modeling, and the results can be extended to more complex settings such as non-ideal gases.  In steady-state, mass flux $\varphi$ is constant, so that equations \eqref{eq:euler:1} will have the form
\begin{align}
  \partial_ x \varphi = 0, \qquad 
  \partial_x p = -\frac{\lambda}{2 D} \frac{\varphi|\varphi|}{\rho}. \label{eq:euler0}
\end{align}
The second equation above defines the change of gas density along the pipe, and integrating it along space yields  
\begin{equation}
  (\rho(L))^2 - (\rho(0))^2 = -\beta\varphi|\varphi|, \qquad \beta = \lambda L / (a^2D). \label{eq:euler0:int}
\end{equation}

\subsection{Pipeline Network Modeling}  \label{sec:network_modeling}

A pipeline network can be modeled as a set of edges $\mathcal{E}$ connected at junctions in a set $\mathcal{V}$ where gas is withdrawn from or injected into the network.  Gas flow is actuated by compressors in a set $\mathcal{C}\subset\mathcal{E}$, which are modeled as multiplicative pressure boosters that preserve through-flow.  Expressing equation \eqref{eq:euler0:int} in terms of squared pressures $\Pi_i$ and $\Pi_j$ at nodes $i,j\in\mathcal{V}$ and flow $\phi_{ij}$ on pipe $(i,j)\in\mathcal{E}$ yields a flow equation for each edge in the set of pipes $\mathcal{P}\subset\mathcal{E}$,
\begin{equation} \label{eq:model1:flowdym} 
  \Pi_j - \Pi_i = \kappa_{ij} \phi_{ij} | \phi_{ij} |, \quad \forall (i,j)\in\mathcal{P},
\end{equation}
where $\kappa_{ij} = a^2 \lambda_{ij} L_{ij} / (A_{ij}^2 D_{ij})$.  The pressure boost achieved by compressor stations is modeled according to
\begin{equation} \label{eq:model1:comps}
 \Pi_j = \alpha_{ij} \Pi_i, \quad \forall (i,j) \in \mathcal{C},
\end{equation}
where $\alpha_{ij}$ is a factor that relates the squares of the compressor discharge and suction pressures.  The network nodes $\mathcal{V}$ are categorized either as slack (pressure) nodes in the set $\mathcal{V}_s$ or nonslack (flow) nodes in the set $\mathcal{V}_q$.  The slack nodes are characterized by a given squared pressure $\Pi_j$, and nonslack nodes are characterized by a given withdrawal flow $q_j$ (negative if an injection), according to
\begin{equation} \label{eq:model1:nodalbc}
    \Pi_j = p_j^2, \,\, j \in \mathcal{V}_s, \quad \text{and} \quad  q_j = d_j-s_j, \,\, j\in \mathcal{V}_q,
\end{equation}
where the values of $p_j$ and withdrawals $d_j\geq 0$ or injections $s_j\geq 0$ must be specified.  We suppose that each node $j\in\mathcal{V}$ is either consumer or supplier, so only one of $d_j$ or $s_j$ can be positive.  Conservation of flow is enforced at each node in the network as
\begin{equation} \label{eq:model1:nodebalance}
  \sum_{i\in\partial_{+}j} \phi_{ij} - \sum_{k\in\partial_{-}j} \phi_{jk} = q_j, \,\, \forall j \in \mathcal{V},
\end{equation}
with the notation $\partial_{+}j=\left\{ i\in \mathcal{V}\mid(i,j)\in \mathcal{E}\right\}$ and $\partial_{-}j=\left\{ k\in \mathcal{V}\mid(j,k)\in \mathcal{E}\right\}$ used to denote the sets of nodes connected to \(j\) by incoming and outgoing edges, respectively.    To ensure mass flow balance for the system, pressure is fixed at a slack node into which flow is free, which represents a large supply source with an injection (negative withdrawal).  The withdrawal $d_j$ or supply $s_j$ at a flow node $j\in\mathcal{V}_q$ specified in \eqref{eq:model1:nodalbc} may be a decision variable.  If the withdrawal $d_j\geq 0$ is optimized, we write that $j\in\mathcal{O}_d\subset\mathcal{V}_q$, and if supply $s_j\geq 0$ is optimized, we write that $j\in\mathcal{O}_s\subset\mathcal{V}_q$.

\subsection{Deterministic Optimal Gas Flow} \label{sec:formulation_det}

The deterministic steady-state optimal gas flow problem is comprised of an objective function solved subject to the  physical flow constraints and boundary conditions defined in equations \eqref{eq:model1:flowdym}-\eqref{eq:model1:nodebalance} as well as inequality constraints, which we define below.  These limitations are imposed due to engineering and operating requirements.  We suppose that the gas pressure at each location $j \in \mathcal{V}$ in the network is bounded according to
\begin{equation} \label{eq:pressure_lim}
    \Pi_j^{min} \leq \Pi_j \leq \Pi_j^{max}, \,\,\, \forall j \in \mathcal{V},
\end{equation}
and that the compressor ratios are bounded according to
\begin{equation} \label{eq:boost_lim}
     1 \leq \alpha_{ij} \leq \alpha_{ij}^{max}, \,\,\, \forall (i,j) \in \mathcal{C}.
\end{equation}
The objective function balances two quantities, the first of which approximates the cost of operating the pipeline by using energy to operate compressors.  Following our previous studies \cite{kazi2024modeling,zlotnik2019optimal}, we approximate the compressor power by the energy needed for adiabatic compression, which takes the form
\begin{equation} \label{eq:obj_comp}
    W_c=\sum_{(i,j)\in\mathcal{C}} \eta_{ij}\phi_{ij}(\alpha_{ij}^m-1), 
\end{equation}
where $0 < m= (\gamma_g-1)/\gamma_g/2 < 1$ where $\gamma_g$ is the heat capacity ratio of the gas \cite{menon05}, and $\eta_{ij}$ is a constant coefficient.  The second component of the objective function quantifies the economic value produced by the pipeline for its users, which is defined as
\begin{equation} \label{eq:obj_econ}
    W_e=\sum_{j\in\mathcal{V}} (c^d_j d_j - c^s_j s_j),
\end{equation}
where $c_j^d$ is the bid price of a consumer and $c_j^s$ is the offer price of a supplier.  In practice, we suppose that the economic value $W_e$ is at least an order of magnitude greater than $W_c$ for an pipeline system design. Finally, in a market mechanism for gas flow scheduling, we may constrain any optimized demand flows by 
\begin{subequations} \label{eq:flow_lims}
\begin{align} 
    0\leq d_j \leq d_j^{\max},  \quad j\in\mathcal{O}_d, \\
    0\leq s_j \leq s_j^{\max},  \quad j\in\mathcal{O}_s.
\end{align}
\end{subequations}
Synthesizing the objective function subject to equality and inequality constraints yields the deterministic steady-state optimal gas flow problem:
\begin{subequations}\label{eq:ogf_det}
\begin{align}
    \min_{\alpha,d,s} & \quad W_c - W_e \,\, \text{as in eq. \eqref{eq:obj_comp} and \eqref{eq:obj_econ}} \label{eq:ogf_det_objective} \\  
    \text{s.t.} & \quad \text{pipe momentum conservation \eqref{eq:model1:flowdym}} \label{eq:ogf_det_momentum_balance} \\    
        & \quad \text{compressor actions \eqref{eq:model1:comps}} \label{eq:ogf_det_compressor_ratio} \\
        & \quad \text{node flow balance \eqref{eq:model1:nodebalance}} \label{eq:ogf_det_flow_balance} \\ 
    & \quad \text{pressure constraints \eqref{eq:pressure_lim}} \label{eq:ogf_det_pressure_bounds}  \\
    & \quad \text{compressor ratio limits \eqref{eq:boost_lim}}  \label{eq:ogf_det_compressor_bounds} \\
    & \quad \text{nomination limits \eqref{eq:flow_lims}}  \label{eq:ogf_det_nom_bounds}
\end{align}  
\end{subequations}
The form of this optimization problem can be modified by collecting variables and parameters into vectors.  Let us define the vectors $\boldsymbol{c}_d$ and $\boldsymbol{c}_s$ of nodal demand and supply prices for optimized flows $\boldsymbol{d}$ and $\boldsymbol{s}$ at nodes in the sets $\mathcal{O}_d$ and $\mathcal{O}_s$, respectively.  Define also the vectors $\boldsymbol{q}$ of nodal gas withdrawals, $\boldsymbol{\kappa}$ of pipe resistances $\kappa_{ij}$, $\boldsymbol{\phi}$ of flows through all pipes, $\boldsymbol{\Pi}$ of all nodal pressures, and $\boldsymbol{\alpha}$ of all compressor ratios. We then define the incidence matrix $\boldsymbol{A}$ of the graph
\begin{equation} \label{eq:incidence}
    A_{ik} = \left\{\begin{array}{ll} 1 & \text{edge $k=(j,i)$ enters node $i$}, \\ 
    -1 & \text{edge $k=(i,j)$ leaves node $i$}, \\ 0 & \text{else}. \end{array}\right.
\end{equation}
The problem \eqref{eq:ogf_det} can then be stated as
\begin{subequations}\label{eq:ogf_det_alt}
\begin{align}
    \min_{\alpha,d,s} & \quad W_c - (\boldsymbol{c}_d^T\boldsymbol{d}-\boldsymbol{c}_s^T\boldsymbol{s}) \label{eq:objective}\\    
    \text{s.t.} & \quad \Pi_j - \Pi_i = \kappa_{ij} \phi_{ij} | \phi_{ij} |, \qquad  \forall (i,j) \in \mathcal{P} \label{eq:det_ogf_alt_momentum_balance}\\ 
    & \quad \Pi_j = \alpha_{ij} \Pi_i, \qquad\qquad\qquad      \forall (i,j) \in \mathcal{C} \label{eq:det_ogf_alt_compressor_ratio}\\
    & \quad \boldsymbol{A} \boldsymbol{\phi} = \boldsymbol{q} \label{eq:det_ogf_alt_flow_balance}\\ 
    & \quad \boldsymbol{\Pi}_{min} \leq \boldsymbol{\Pi} \leq \boldsymbol{\Pi}_{max}  \label{eq:det_ogf_alt_pressure_bounds} \\
    & \quad \boldsymbol{1} \leq \boldsymbol{\alpha} \leq \boldsymbol{\alpha}_{max}   \label{eq:det_ogf_alt_compressor_bounds} \\
    & \quad \boldsymbol{0} \leq \boldsymbol{d} \leq \boldsymbol{d}_{max}, \,\, \boldsymbol{0} \leq \boldsymbol{s} \leq \boldsymbol{s}_{max}     \label{eq:det_ogf_alt_nom_bounds}
\end{align}  
\end{subequations}
Formulations similar to the optimization problem \eqref{eq:ogf_det_alt} have been examined in previous studies \cite{vuffray2015monotonicity,wu2017adaptive}, and extended to transient flows \cite{zlotnik2019optimal} and gas mixtures \cite{kazi2024modeling}, and economic interpretations were also developed \cite{rudkevich17hicss,sodwatana2023optimization}.  Steady-state formulations such as problem \eqref{eq:ogf_det} can be viewed as providing a day-ahead flow allocation for an operating day, assuming steady ratable offtakes by customers.  Here we extend the steady-state optimal gas flow problem to account for uncertain demands using the SFV approach.

\section{Chance-Constrained Optimal Gas Flow} \label{sec:ccssogf}

In the steady-state flow setting, we aim to develop a chance-constrained formulation to provide probabilistic guarantees given a known distribution in one or more gas flow withdrawals from a pipeline network.  This type of formulation can be more flexibly calibrated than robust optimization approaches that compensate for arbitrary interval uncertainty \cite{vuffray2015monotonicity,misra2020monotonicity}. The key concern for pipeline systems is to assure delivery of gas to all customers given uncertainty in some loads while maintaining adequate pipeline pressure.  The lower bound in constraint \eqref{eq:pressure_lim} (i.e., \eqref{eq:det_ogf_alt_pressure_bounds}) is somewhat flexible in practice, as long as it is not violated significantly and pressures can be returned to adequate levels for the subsequent operating day.  We enforce this limit using chance constraints, allowing for a small probability $\epsilon$ of violation. 

Suppose that a subset of gas pipeline nodes $\mathcal{S}\subset \mathcal{V}$ has stochastic gas consumptions of the form
\begin{equation} \label{eq:stochastic_load}
    d_j(\omega) = d_j + r_j(\omega), \quad j\in\mathcal{S},
\end{equation}
where we suppose that $(\Omega_j,\mathcal{B}_j,\mu_j)$ is a probability space where $r_j:\Omega_j\to R_j$ is a random variable taking values on a compact interval $R_j=[\underline{r}_j,\overline{r}_j]\subset\mathbb{R}$ and where $\mathcal{B}_j$ is the Borel $\sigma$-algebra.  The fixed parameters $d_j$ denote the nominal baseline load.  Letting $\Omega=\prod_{j=1}^{|\mathcal{S}|}\Omega_j$, we consider samples $\omega\in\Omega$ as instances of stochastic loads across the pipeline network.  We suppose that values of random parameters $r_j(\omega)$ corresponding to various values of $\boldsymbol{q}$ in problem \eqref{eq:ogf_det_alt} would result in different optimal solutions of that problem.  Thus in the stochastic extension below we consider the physical solutions $\boldsymbol{\Pi}$ and $\boldsymbol{\phi}$ to also depend on the sample $\omega$, and we seek to find a single decision for the value of $\boldsymbol{\alpha}$ and any optimized flows $d_j$ or $s_j$ for $j\in\mathcal{O}\subset\mathcal{V}$ such that the inequality constraints corresponding to \eqref{eq:det_ogf_alt_pressure_bounds} are satisfied in a probabilistic sense.    The flow equations \eqref{eq:model1:flowdym} and compressor conditions \eqref{eq:model1:comps} must be evaluated for all samples $\omega$ as
\begin{align} \label{eq:model1:flowdym_rand} 
 \!\!\!\!\!\!\Pi_j(\omega) \!-\! \Pi_i(\omega) \!=\! \kappa_{ij} \phi_{ij}(\omega) | \phi_{ij}(\omega) |, \, \forall (i,j)\!\in\!\mathcal{P}, \forall \omega\!\in\!\Omega, \!\!\!\! &
\\ \label{eq:model1:comps_rand}
 \Pi_j(\omega) = \alpha_{ij} \Pi_i(\omega), \quad \forall (i,j) \in \mathcal{C}, \,\, \forall \omega\in\Omega. &
\end{align}
Conservation of flow \eqref{eq:model1:nodebalance} is enforced for all samples as
\begin{equation} \label{eq:model1:nodebalance_rand}
  \sum_{i\in\partial_{+}j} \phi_{ij}(\omega) \!-\!\!\! \sum_{k\in\partial_{-}j} \phi_{jk}(\omega) \!=\! q_j(\omega), \,\, \forall j \!\in\! \mathcal{V}, \,\, \forall \omega\!\in\!\Omega,
\end{equation}
where $q_j(\omega) = d_j+r_j(\omega)-s_j$.  The chance constraint for minimum pipeline pressure is expressed using a quadratic penalty function 
\begin{subequations} \label{eq:penalty}
  \begin{equation}
    \Gamma(z) = \begin{cases}
        \gamma z^2, \quad& \mbox{if } x \geq 0, \\
        0, \quad & \mbox{otherwise},
    \end{cases}
\end{equation}
and a penalized minimum pressure violation variable
\begin{equation}\label{eq:pressure_violation}
    v^{min}_j(\omega) = \Gamma\left(\Pi_j^{\min} - \Pi_j(\omega)\right).
\end{equation}
\end{subequations}
The quadratic penalty is chosen in order to reflect the practical consideration that larger violations of the lower pressure limits are more problematic. Additionally, because it is twice differentiable, it facilitates a well-behaved nonlinear program.  The chance constraint is then
\begin{equation}\label{eq:chance_constraint}
    \mathbb{E_{\omega}} \!\left[v^{min}_j(t,\omega) \right] \leq \epsilon_j, \qquad \forall t \in T,
\end{equation}
where the acceptable violation probability $\epsilon_j$ may depend on the network node $j\in\mathcal{V}$.
We approximate the constraint \eqref{eq:chance_constraint} in a deterministic manner (i.e., without requiring Monte Carlo simulation) by applying the SFV method as follows.  We discretize the stochastic space $\Omega$, which has one-to-one correspondence with the interval $R_j$ on which the random consumption $r_j(\omega)$ appears, into $M$ cells delimited by $M+1$ uniformly spaced boundary points, each of which corresponds to a value of $r_j{\omega_m}$.  We construct the penalty variable using a third order spline expansion on $\Omega$ of form
\begin{subequations}
\begin{equation}\label{eq:violation_interpolation}
    v^{min}_j(\omega) =  \Gamma\left(\Pi_j^{\min} - \Pi_j(\omega)\right) = \sum_{m \in M} a_{jm} b_{jm}(\omega),
\end{equation}
where $b_{jm}(\omega)$ is the $m$-th spline function on the uniform stochastic space grid that is completely known.  The constraint \eqref{eq:chance_constraint} is expressed as an expectation over $\omega$ as
\begin{equation}\label{eq:b_spline_chance_constraint}
    \mathbb{E}_\omega[v^{min}_j(\omega)]=\sum_{m \in M} a_{jm} \int_{\omega} b_{jm}(\omega) d \mu_j \leq \epsilon_j,
\end{equation}
\end{subequations}
where the coefficients $a_{jm}$ serve as decision variables that can be optimized.  Aside from minimum pressure limits, all other constraints including the maximum operating pressure,
\begin{equation} \label{eq:pressure_lim_up_rand}
    \Pi_j(\omega) \leq \Pi_j^{max}, \,\,\, \forall j \in \mathcal{V},
\end{equation}
are strictly enforced for all values of $\omega$.    The objective function components must also be evaluated in expectation,
\begin{align} \label{eq:obj_comp_exp}
    \mathbb{E}_{\omega}[W_c] & =\sum_{(i,j)\in\mathcal{C}} \eta_{ij}(\alpha_{ij}^m-1) \cdot \mathbb{E}_{\omega}[\phi_{ij}], \\
\label{eq:obj_econ_exp}    \mathbb{E}_{\omega}[W_e] & =\sum_{j\in\mathcal{V}} (c^d_j (d_j + \mathbb{E}_{\omega}[r_j]) - c^s_j s_j).
\end{align}
Assembling the objective function and probabilistic constraints in the stochastic setting results in the chance-constrained optimal gas flow problem given by
\begin{subequations}\label{eq:ogf_rand}
\begin{align}
    \min_{\alpha,d,s} & \quad \mathbb{E}_{\omega}[W_c] - \mathbb{E}_{\omega}[W_c] \,\, \text{as in eq. \eqref{eq:obj_comp_exp} and \eqref{eq:obj_econ_exp}} \label{eq:ogf_rand_objective_exp} \\    
    \text{s.t.} & \quad \text{pipe momentum conservation \eqref{eq:model1:flowdym_rand} } \label{eq:ogf_rand_momentum_balance} \\    
        & \quad \text{node flow balance \eqref{eq:model1:nodebalance_rand}} \label{eq:ogf_rand_flow_balance} \\ 
    & \quad \text{compressor actions \eqref{eq:model1:comps_rand}} \label{eq:ogf_rand_compressor_ratio} \\
    & \quad \text{pressure constraints \eqref{eq:violation_interpolation}-\eqref{eq:b_spline_chance_constraint} and \eqref{eq:pressure_lim_up_rand}} \label{eq:ogf_rand_pressure_bounds}  \\
    & \quad \text{compressor ratio limits \eqref{eq:boost_lim}}  \label{eq:ogf_rand_compressor_bounds} \\
    & \quad \text{nomination limits \eqref{eq:flow_lims}}  \label{eq:ogf_rand_nom_bounds}
\end{align}  
\end{subequations}
Invoking the vector notations used to express formulation \eqref{eq:ogf_det} as \eqref{eq:ogf_det_alt} above, and using $\boldsymbol{r}$ to state the vector of random withdrawals provided at prices $\boldsymbol{c}_r$ for $j\in\mathcal{S}$, we re-write problem \eqref{eq:ogf_rand} as
\begin{subequations}\label{eq:ogf_rand_alt}
\begin{align}
    \!\!\!\!\!\!\!\! \min_{\alpha,d,s} & \,\, \mathbb{E}_{\omega}[W_c] - (\boldsymbol{c}_d^T\boldsymbol{d}-\boldsymbol{c}_s^T\boldsymbol{s}) - \boldsymbol{c}_r^T\mathbb{E}_{\omega}[\boldsymbol{r}] \label{eq:objective_rand_alt}\\    
    \text{s.t.} & \,\,  \Pi_j(\omega) \! -\!  \Pi_i(\omega) \! =\!  \kappa_{ij} \phi_{ij}(\omega) | \phi_{ij}(\omega) |,  \nonumber \\ & \qquad\qquad\qquad\qquad\qquad \,\,  \forall (i,j) \in \mathcal{P}, \,\, \forall \omega\in\Omega, \label{eq:rand_ogf_alt_momentum_balance}\\ 
    & \,\, \Pi_j(\omega) = \alpha_{ij} \Pi_i(\omega), \quad\quad  \forall (i,j) \in \mathcal{C}, \,\, \forall \omega\in\Omega, \label{eq:rand_ogf_alt_compressor_ratio}\\
    & \,\, \boldsymbol{A} \boldsymbol{\phi}(\omega) = \boldsymbol{q}(\omega), \quad \forall \omega\in\Omega \label{eq:rand_ogf_alt_flow_balance}\\ 
        & \,\, \boldsymbol{\Pi}(\omega) \leq \boldsymbol{\Pi}_{max}, \quad \forall \omega\in\Omega    \label{eq:rand_ogf_alt_pressure_bounds_up} \\
    & \,\, \boldsymbol{1} \leq \boldsymbol{\alpha} \leq \boldsymbol{\alpha}_{max} \label{eq:rand_ogf_alt_compressor_bounds} \\
    & \,\, \Gamma\left(\Pi_j^{\min} - \Pi_j(\omega)\right) = \sum_{m \in M} a_{jm} b_{jm}(\omega), \nonumber \\ & \qquad\qquad\qquad\qquad\qquad\qquad \,\, \forall j\in\mathcal{V}, \, \forall \omega\in\Omega,  \,\, \label{eq:rand_ogf_alt_viol_interp} \\
    & \sum_{m \in M}\!\! a_{jm} \int_{\omega} b_{jm}(\omega) d \mu_j \leq \epsilon_j,
    \,\, \forall j\!\in\!\mathcal{V}, \, \forall \omega\!\in\!\Omega,  \,\,  \,\,  \,\,  \label{eq:rand_ogf_alt_b_spline_cc} \\
    & \,\, \boldsymbol{0} \leq \boldsymbol{d} \leq \boldsymbol{d}_{max}, \,\, \boldsymbol{0} \leq \boldsymbol{s} \leq \boldsymbol{s}_{max}     \label{eq:rand_ogf_alt_nom_bounds}
\end{align}  
\end{subequations}
In the rest of the paper, we examine the computational solution of problem \eqref{eq:ogf_rand_alt} for to small test networks by estimating the distributions of physical solutions.  Because the optimization objective \eqref{eq:objective_rand_alt} is formulated as an expected economic value, we also examine the possibility of reconstructing the distributions of dual variables and interpreting them as stochastic locational values of natural gas.

\vspace{-0.5ex}
\section{Computational Studies} \label{sec:examples}
\vspace{-0.5ex}

We explore solutions to problem \eqref{eq:ogf_rand_alt} obtained using our SFV-based uncertainty management method using two example systems.  The first consists of a single pipe with a compressor at the sending end, and the second is an 8-node test network that consists of 5 pipes and 3 compressor stations, and which was described in our previous studies \cite{gyrya2019explicit,tokareva2024stochastic,kazi2024modeling}. We examine the distributions of physical flows and pressures that arise given uniform or truncated normal probability distributions in withdrawal gas flows. The problem \eqref{eq:ogf_rand_alt} is implemented as a nonlinear constrained optimization model in JuMP \cite{dunning2017jump}, an algebraic modeling language in Julia, and solving using an open source interior point method based optimization solver IPOPT \cite{waechter2006ipopt}.  For computational well-conditioning, we re-scale the pressure, flow, distance, and velocity variables according to the non-dimensionalization in a previous study \cite{srinivasan2022numerical}.   We suppose for simplicity in this study that the value of the exponent $m$ in equation \eqref{eq:obj_comp} is $m = 1$, so that the specific heat capacity ratio of the gas is 0.5.
\vspace{-1.5ex}
\begin{figure}[!h]
    \centering
    \includegraphics[width=\linewidth]{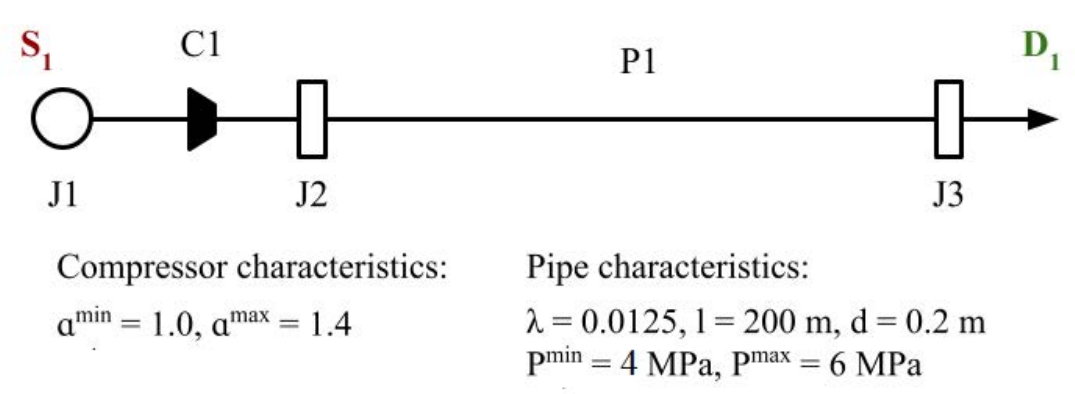}
    \vspace{-5ex}
    \caption{Single pipe test system and component parameters.}
    \label{fig:3-node}
\end{figure}
\vspace{-2ex}

\subsection{Single Pipe with Compressor} \label{sec:singlepipe}

\begin{figure}[!h]
    \centering
    \includegraphics[width=\linewidth]{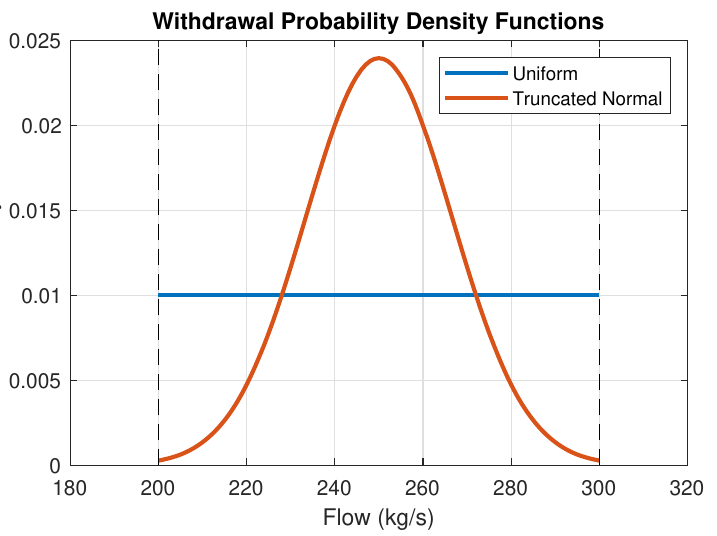}
        \vspace{-5ex}
    \caption{Single pipe example probability distribution functions $q_3\sim \mathcal{U}[200,300]$ and $q_3\sim N(200,50/3)$ with truncated tails.}
    \label{fig:withdrawal_pdfs}
\end{figure}

\begin{figure}[!t]
    \centering
    \includegraphics[width=\linewidth]{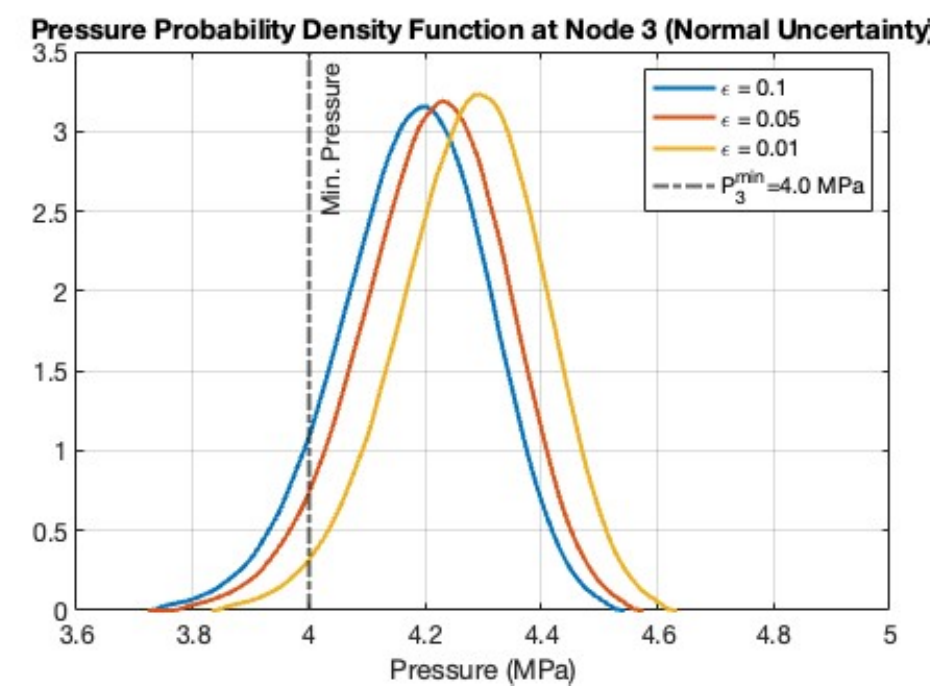}
    \includegraphics[width=\linewidth]{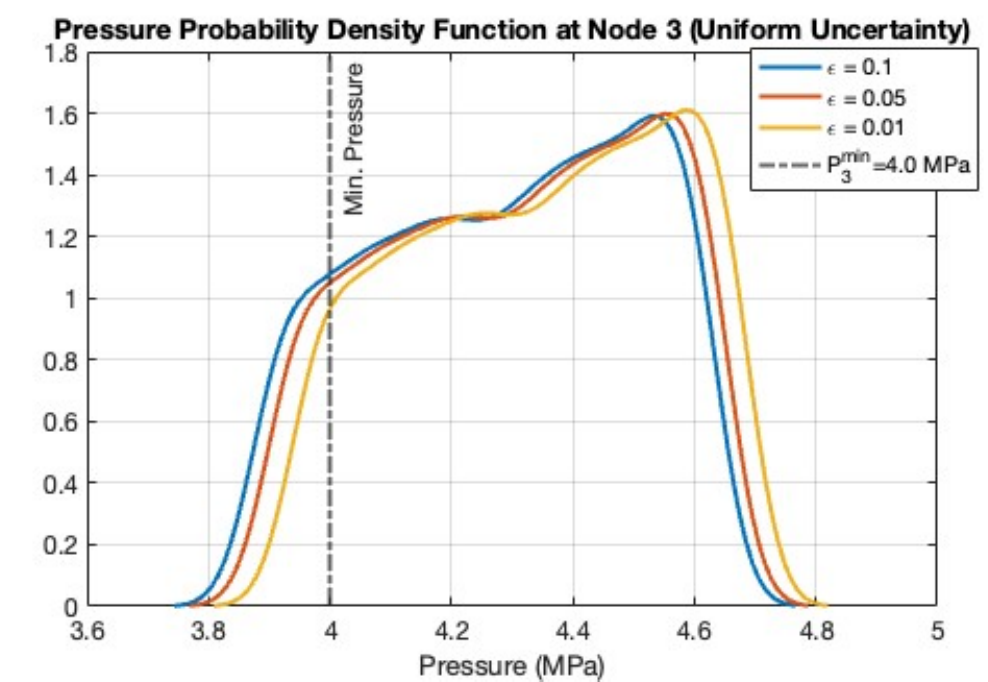}
                  \vspace{-4ex}
    \caption{Pressure Probability Distribution at Withdrawal Node $J3$. Top: Truncated Normal Uncertainty $q_3\sim N(200,50/3)$;  Bottom: Uniform Uncertainty $q_3\sim \mathcal{U}[200,300]$}
    \label{fig:pressure_distributions_singlepipe}
\end{figure}

The first test network consists of a single pipe with a compressor connected to the inlet node, as shown in Figure \ref{fig:3-node}. The injection node (N1) acts as a slack node with fixed pressure set at $p^{slack}_1 = 4.3367$ MPa. The pressure bounds at the withdrawal node (N3) are specified as $P_{min} = 4.0$ MPa and $P_{max} = 6.0$ MPa. The compressor ratio is bounded by $\alpha \in [1,1.4]$. The uncertain withdrawal flow is parameterized by an uncertain parameter ($\omega$) as $q_3(\omega) = d_3+r_3(\omega)$, where $d_3$ represents the nominal value that is fixed at 250 kg/s.  We consider a scenario where the uncertain flow $r_3(\omega)$ follows a uniform distribution $r_3\sim \mathcal{U}[-50,50]$, and a scenario where $r_3(\omega)$ follows a truncated normal distribution with zero mean and standard deviation $\sigma=50/3$, i.e., $r_3\sim \mathcal{N}(0,50/3)$ with support on $[-50,50]$.  The distributions of the withdrawal $q_3(\omega)$  are shown in Figure \ref{fig:withdrawal_pdfs}.
\vspace{-1ex}
\vspace{-1ex}
We discretize the sample space $\Omega_3$ into $K=100$ stochastic cells that uniformly partition the interval $[200,300]$ and solve the reformulated chance constrained optimization problem \eqref{eq:ogf_rand_alt} for several values of the acceptable violation probability parameter $\epsilon_3$, namely, $\epsilon_3 \in\{0.01,0.05,0.1\}$. The distribution of the pressure at the withdrawal node N3 is approximated by the kernel density estimate method (\verb"ksdensity" in MATLAB) using 10000 empirical samples of pressure distribution obtained from the optimal solution. The pressure probability density function (PDF) estimate at the withdrawal node is plotted for these constraint violation probabilities and shown in Figure \ref{fig:pressure_distributions_singlepipe}. 

Observe that the withdrawal pressure PDF retains a bell curve type shape for normally distributed uncertain withdrawal flow but the pressure PDF for uniform uncertain withdrawal has a non-uniform distribution over the pressure values. The withdrawal pressure PDFs are consistent as they show higher violation for higher values of $\epsilon$ and vice-versa. We also note that the estimated probability of violation for the pressure constraint ($P \leq P_{min}$) is higher in the uniform uncertainty case than the normal distribution for same chance constraint violation probability parameter ($\epsilon_3$). This happens because the pressure drop inside the pipe increases with increasing withdrawal flow resulting in lower pressure at the withdrawal node, and the probability of withdrawal flow near the upper end of the uncertainty interval $[200,300]$ is higher given uniform uncertainty than when considering the normally distributed uncertainty scenario. A uniform distribution for the uncertain flow also results in higher compressor ratio solution $\alpha = \{1.1985, 1.192, 1.187\}$ compared to normal uncertainty ($\alpha = \{1.183, 1.171, 1.164\}$), where these enumerations correspond to the values $\epsilon_3=\{0.01,0.05,0.1\}$ for the acceptable violation probability.

\vspace{-1.5ex}
\begin{figure}[!h]
    \centering
    \includegraphics[width=\linewidth]{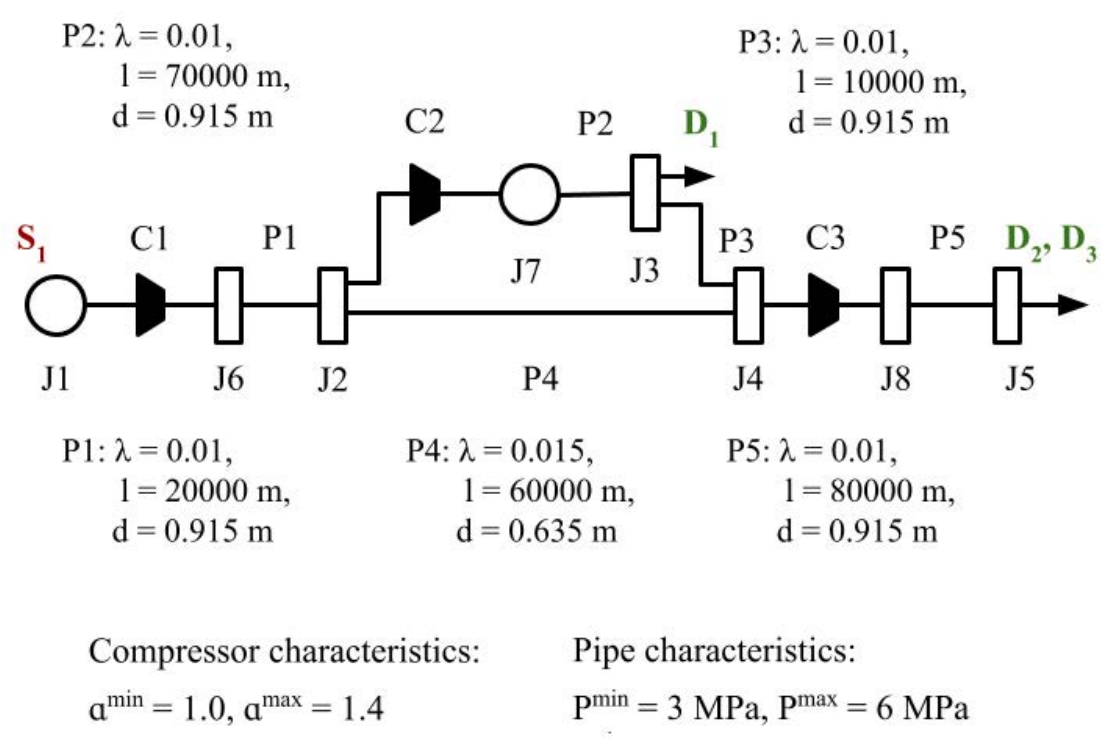}
    \vspace{-4ex}
    \caption{8-node Gas Pipeline Network System}
    \label{fig:8-node}
\end{figure}
\vspace{-3ex}

\subsection{8-node Pipeline Test Network} \label{sec:eightnode}

In our second case study, we apply the SFV-based optimization approach to an uncertainty management scenario for an 8-node gas pipeline test network with 5 pipes, 3 compressors, 2 withdrawal (at nodes J3 and J5), and 1 slack injection node (J1) with constant pressure $P^{slack}_1 = 5.0$ MPa (shown in Fig.\ref{fig:8-node}).  This test network has been examined for steady-state and transient flow algorithms in several previous studies \cite{gyrya2019explicit,tokareva2024stochastic, sodwatana2023optimization,kazi2024modeling}. The lower bounds for pressure at withdrawal nodes is specified as $P^{min}_5 = 4.0$ MPa at node J5, and $P^{min}_3 = 3.0$ MPa at the other nodes. We suppose that the withdrawal flow at node J5 is uncertain and parameterized using a uniformly distributed uncertain parameter, whereas the withdrawal flow $q_{3}$ at node J3 is a decision variable in the optimization problem. The uncertain withdrawal flow at node J5 is modeled as $q_5(\omega)=d_5+r_5(\omega)$, where $d_5=64$ kg/s and $r_j = \mathcal{U}[0,32]$.  In our implementation, $\Omega_5=[0,32]$ is discretized using $K=50$ uniform stochastic cells. We solve three cases with different values of the withdrawal nomination constraint bound $q_3^{\max}\equiv d_3^{\max}$, namely, $q_3^{\max}\equiv d^{max}_{3}\in \{200,300,\infty\}$ kg/s, using the same value of violation probability parameter $\epsilon = 0.1$ in each case. 

\begin{figure}
    \centering
    \vspace{1ex}
    \includegraphics[width=\linewidth]{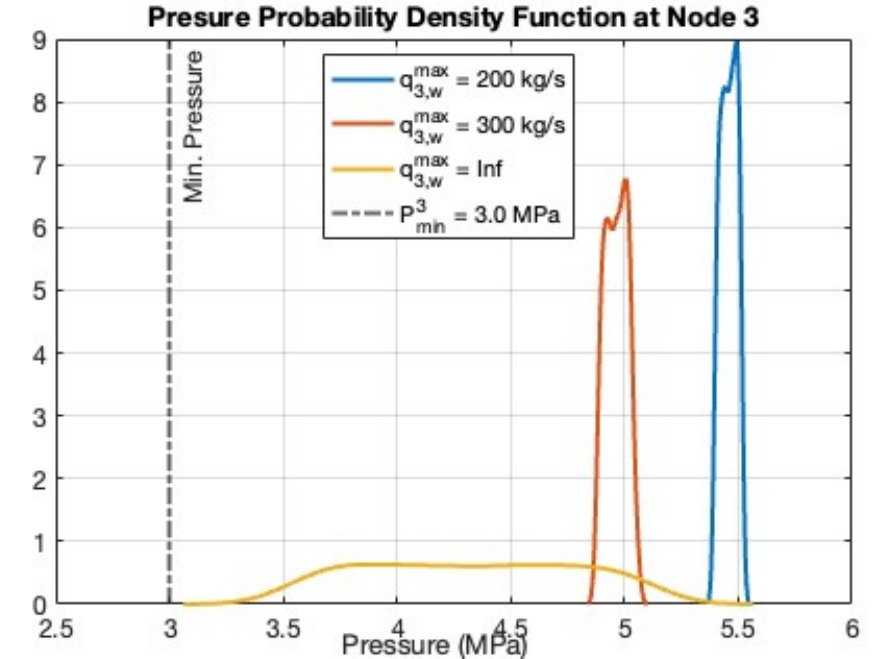}
    \vspace{1ex}
    \includegraphics[width=\linewidth]{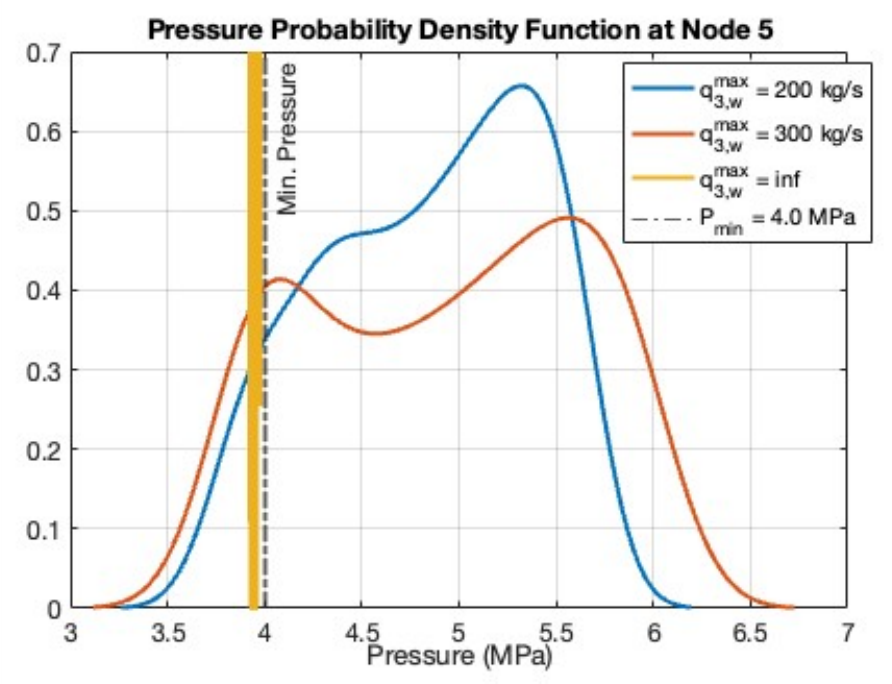}
        \vspace{-5ex}
    \caption{Pressure probability distributions for 8-node case study. Top: Node 3; and Bottom: Node 5.}
    \label{fig:pressure_nodel_8node}
\end{figure}

\begin{figure}[!t]
    \centering
    \vspace{-2ex}
    \includegraphics[width=\linewidth]{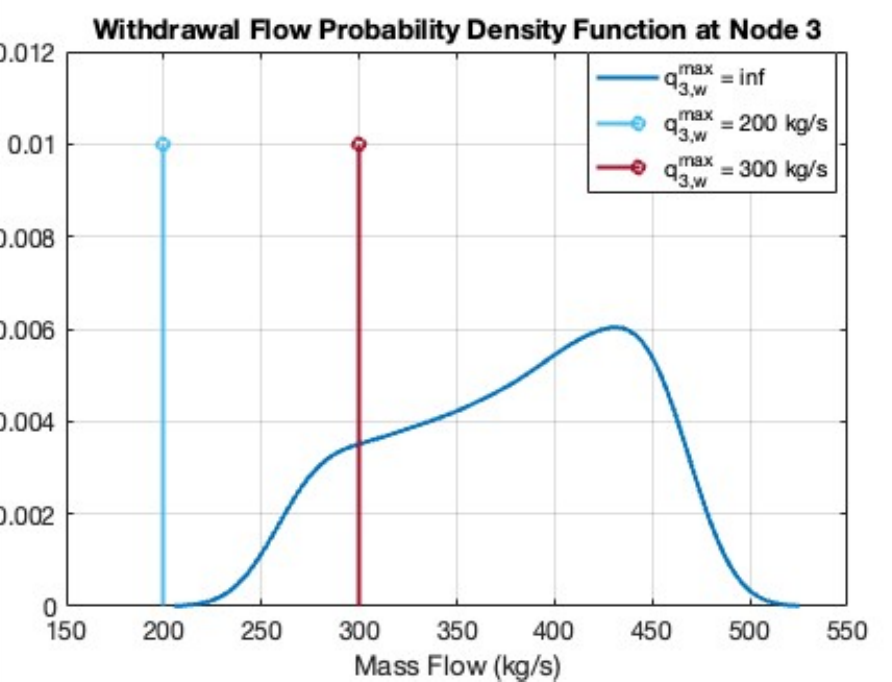}
    \vspace{-4ex}
    \caption{Withdrawal Flow Probability Distribution at Node 3}
    \label{fig:withdrawal_node3}
\end{figure}

Estimates of the pressure PDFs at two withdrawal nodes J3 and J5 are shown in Figure \ref{fig:pressure_nodel_8node}. The pressure profile shows no violation at withdrawal node J3 for any instance. In contrast, the pressure at withdrawal node J5 exhibits bound violation for all instances and the pressure profile for the instance with no upper bound ($q_3^{max} = \infty$) has values less than $P^{min}_5 = 4.0$ MPa for each uncertain scenario. The pressure PDFs also show a higher standard deviation in withdrawal pressure at node J3 for $q^{max}_{3} = \infty$ and in withdrawal pressure at node J5 for $q_3^{max} \in \{200,300\}$ kg/s.  
The optimal withdrawal flow $q_{3}$ at J3 is shown in Fig \ref{fig:withdrawal_node3}. We observe that the withdrawal is at upper bound limits for $q_3^{\max}\equiv d_3^{\max} \in \{200,300\}$ kg/s, but has a continuous PDF estimate when there is no upper bound for the withdrawal flow variable. 

\subsection{Economic Interpretation of Lagrange Multipliers}  \label{sec:economic}

\begin{figure}
    \centering
    \includegraphics[width=\linewidth]{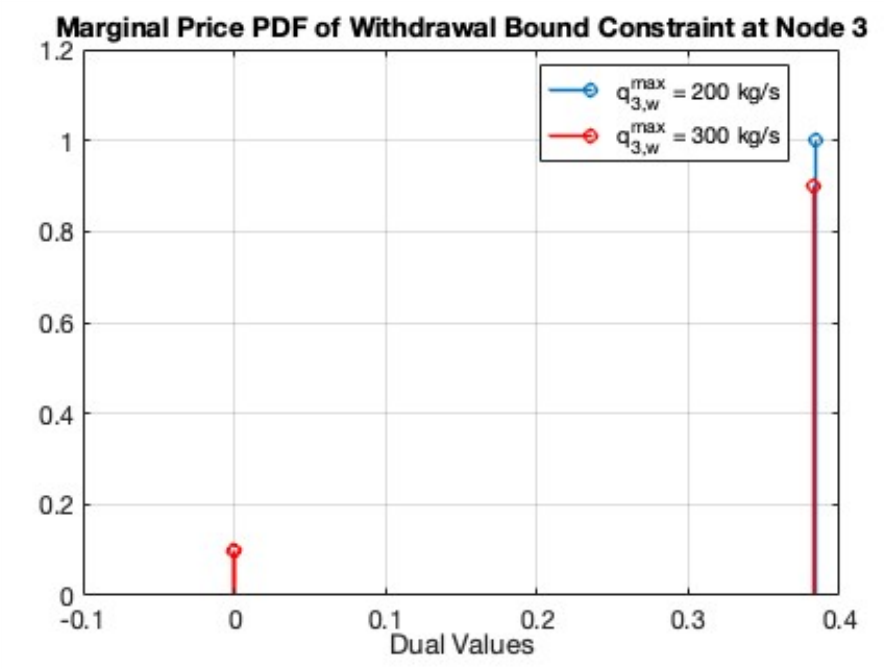}
    \includegraphics[width=\linewidth]{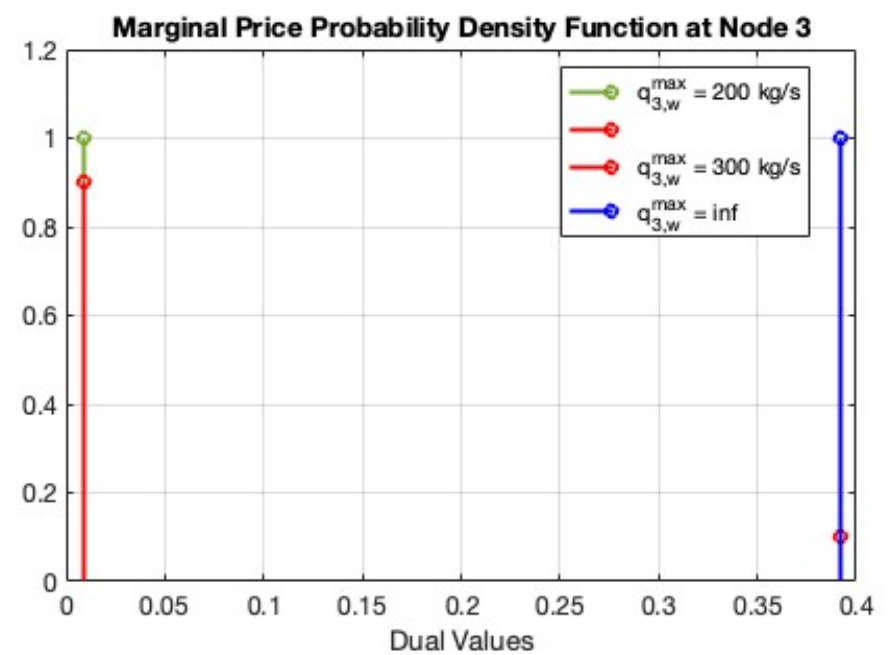}
    \caption{Marginal price probability distributions at Node 3. Top: associated with the withdrawal bound constraint; Bottom: associated with the flow balance constraint}
    \label{fig:duals_node3}
\end{figure}

We now explore the possibility of economic analysis of solutions to problem \eqref{eq:ogf_rand_alt} by use of the associated Lagrange multipliers to reconstruct stochastic locational prices.  In particular, by solving problem \eqref{eq:ogf_rand_alt} we obtain the Lagrange multipliers $\boldsymbol{\lambda}_q(\omega)$ that are associated with the nodal mass flow balance constraint \eqref{eq:rand_ogf_alt_flow_balance}, as well as Lagrange multipliers $\boldsymbol{\lambda}_d(\omega)$ associated with upper bound constraint \eqref{eq:rand_ogf_alt_nom_bounds} for the optimized withdrawal flows the withdrawal bound constraint (in particular, $d_{3}(\omega) \leq d^{max}_{3}$ the 8-node example). These multipliers at the optimal solution, which we interpret as marginal prices, follow the Karush-Kuhn Tucker (KKT) conditions, and can be derived by assembling terms in the Lagrangian as shown in equation \eqref{eq:stoch_lagrangian}, differentiating with respect to $\boldsymbol{d}$ to obtain the first order condition in equation \eqref{eq:stoch_kkt_firstorder}, and solving for the extremal value. For the 8-node example, the relevant parts of the Lagrangian take the form
\begin{align} 
        &L(.,\boldsymbol{d}) = - c_3^d(d_3+\mathbb{E}_{\omega}[r_3]) + \sum_{\omega \in \Omega} \boldsymbol{\lambda}_q^T(\boldsymbol{A} \boldsymbol{\phi}(\omega) - \boldsymbol{q}(\omega)) \nonumber \\ & \qquad \qquad +...+ \sum_{\omega \in \Omega} \lambda_{d,3} (d_{3} - d^{max}_{3}), \label{eq:stoch_lagrangian}
\end{align}
and taking the variation with respect to the optimized withdrawal $d_3$ yeilds
\begin{align}
        &  \quad \nabla L_{q_3}(.) = - \frac{c_3}{K} + \lambda_{q,3} + \lambda_{d,3} = 0. \label{eq:stoch_kkt_firstorder}
\end{align}
Recall that the parameter $K$ denotes the number of stochastic cells used to discretize the interval $\Omega$, which in this case is $\Omega_3\equiv[200,300]$.
\begin{figure}[!t]
    \centering
    \includegraphics[width=\linewidth]{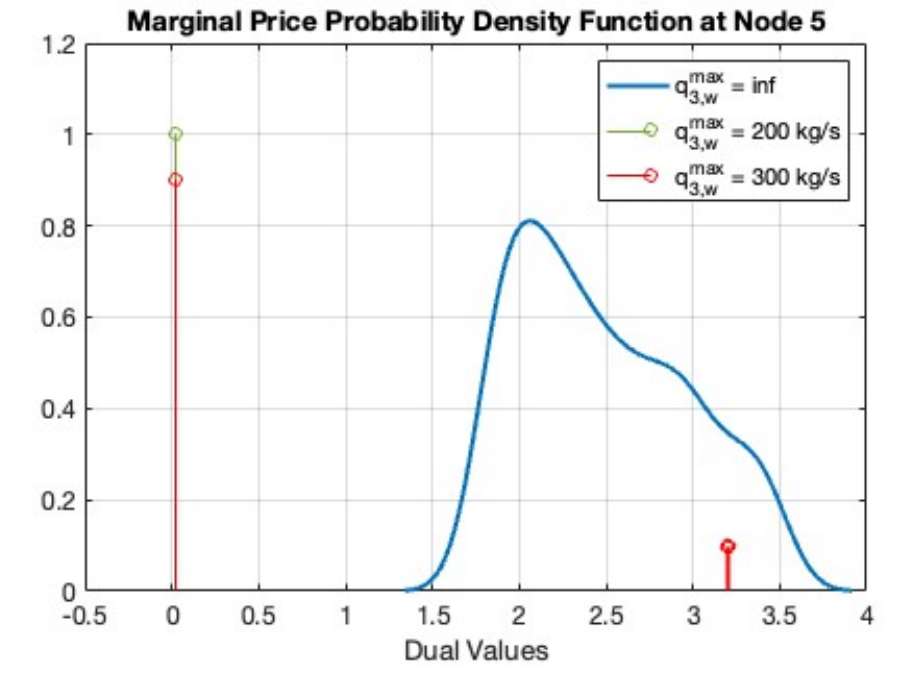}
    \caption{Marginal price probability distribution associated with the flow balance constraint at Node J5}
    \label{fig:dual_flow_node5}
\end{figure}
The probability distribution functions of $\lambda_{f,3}$ and $\lambda_{q,3}$ are shown in Figure \ref{fig:duals_node3}, and in this case they are discrete.  Observe that the dual variables satisfy the KKT condition where $c_3 = 20$ and $K=50$ and the sum of the dual variables is equal to $c_3/K = 0.4$. We also observe that for instances with $q^{max}_{3}\in\{200,\infty\}$ kg/s, the dual variable $\lambda_{q,3}$ is constant at 0.02 and 0.39 respectively for each scenario, whereas the dual variable $\lambda_{d,3} \approx$ 0 and 0.38 respectively. For the instance with $q^{max}_{3}=300$ kg/s, the dual variables ($\lambda_{q,3},\lambda_{d,3}$) are dependent on the scenario, and their sum is equal to 0.4  as expected based on \eqref{eq:stoch_kkt_firstorder}.  Based on these observations, we see that the SFV-based computational method satisfies the KKT conditions for this case study.

Let us now examine the PDF estimate for locational marginal price $\lambda_{q,5}$ at withdrawal node J5, which is shown in Figure \ref{fig:dual_flow_node5}. The marginal price denotes the sensitivity of the optimal objective function value $J(\alpha,d,s)$ to variation in the uncertain withdrawal flow $q_{5}(\omega)$ as $\lambda_{q,5}(\omega) = \partial J(\alpha,d,s)/\partial q_{5}(\omega)$. Because $q_{5}(\omega)$ is a random variable, then we the multiplier $\lambda_{q,5}(\omega)$ is a random variable as well, and its distribution can be reconstructed from the dual variables corresponding to flow balance constraints \eqref{eq:rand_ogf_alt_flow_balance} in the optimal solution to problem \eqref{eq:ogf_rand_alt}.   We observe that distribution of the marginal price $\lambda_{q,5}(\omega)$ of gas at node J5 is always zero for $q_3^{\max}=200$ and sometimes 3.18 when $q_3^{\max}=300$.  If there is no specified upper bound on optimized withdrawal to node J3, i.e., $q_3^{\max}=\infty$, then the optimization solution is driven by the probabilistic lower pressure consraint at node J5  and as a result the marginal price distribution at node J3 is continuous.

\section{Conclusion} \label{sec:conc}

We have developed a stochastic finite volume representation for solutions of the steady-state gas flow equations on a pipeline network with compressors subject to uncertain boundary conditions.  This representation is used to formulate a stochastic optimization problem for steady-state flow allocation in a pipeline system subject to uncertain demands. This type of uncertainty modeling enables the generalizable representation of arbitrary probability distributions in problem parameters, and facilitates higher-order reconstructions in computationally efficient calibration of probabilistic inequality constraints.  The primary contribution of our study is a chance-constrained optimization formulation for steady-state flow allocation for a gas pipeline system, and a method for computing solutions for compressor configuration and nomination flows that optimize an operational and/or economic objective in expectation.  By extending the proven method of chance-constraints using higher-order nonlinear reconstructions, our pipeline demand uncertainty management approach provides probabilistic robustness guarantees that can be calibrated based on risk tolerance.  Furthermore, when the objective function is interpreted as an expected economic value created for users of the pipeline system, we show that the dual variables of the optimization solution can be used to reconstruct distributions that we interpret as probability densities of locational gas prices.  Such solutions could be used to quantify the financial volatility of energy in an optimization-driven gas market mechanism that accounts for demand, supply, transportation capacity, and price volatility.  Compelling future work would include analysis of convergence in modeling accuracy and probability measure, as well as an investigation of the trade-offs between computational cost and stochastic degrees of freedom examined.  We note that, because consumption uncertainty throughout a network may be correlated, it is possible to represent wide-area uncertainty using only a few stochastic dimensions.

\bibliographystyle{unsrt}
\bibliography{references.bib}

\begin{thebibliography}{10}

\bibitem{yang2017effect}
Jingwei Yang, Ning Zhang, Chongqing Kang, and Qing Xia.
\newblock Effect of natural gas flow dynamics in robust generation scheduling
  under wind uncertainty.
\newblock {\em IEEE Transactions on Power Systems}, 33(2):2087--2097, 2017.

\bibitem{roald2017chance}
Line Roald and G{\"o}ran Andersson.
\newblock Chance-constrained {AC} optimal power flow: Reformulations and
  efficient algorithms.
\newblock {\em IEEE Transactions on Power Systems}, 33(3):2906--2918, 2017.

\bibitem{zhao2016unit}
Bining Zhao, Antonio~J Conejo, and Ramteen Sioshansi.
\newblock Unit commitment under gas-supply uncertainty and gas-price
  variability.
\newblock {\em IEEE Transactions on Power Systems}, 32(3):2394--2405, 2016.

\bibitem{levary1980natural}
Reuven~R. Levary and Burton~V. Dean.
\newblock A natural gas flow model under uncertainty in demand.
\newblock {\em Operations Research}, 28(6):1360--1374, 1980.

\bibitem{sahinidis2004opportunities}
Nikolaos~V. Sahinidis.
\newblock Optimization under uncertainty: state-of-the-art and opportunities.
\newblock {\em Computers \& chemical engineering}, 28(6-7):971--983, 2004.

\bibitem{bertsimas2011theory}
Dimitris Bertsimas, David~B. Brown, and Constantine Caramanis.
\newblock Theory and applications of robust optimization.
\newblock {\em SIAM review}, 53(3):464--501, 2011.

\bibitem{behrooz2016managing}
Hesam~Ahmadian Behrooz.
\newblock Managing demand uncertainty in natural gas transmission networks.
\newblock {\em Journal of Natural Gas Science and Engineering}, 34:100--111,
  2016.

\bibitem{zavala2014stochastic}
Victor~M. Zavala.
\newblock Stochastic optimal control model for natural gas networks.
\newblock {\em Computers \& Chemical Engineering}, 64:103--113, 2014.

\bibitem{liu2020dynamic}
Kai Liu, Lorenz~T. Biegler, Bingjian Zhang, and Qinglin Chen.
\newblock Dynamic optimization of natural gas pipeline networks with demand and
  composition uncertainty.
\newblock {\em Chemical Engineering Science}, 215:115449, 2020.

\bibitem{roald2020uncertainty}
Line~A. Roald, Kaarthik Sundar, Anatoly Zlotnik, Sidhant Misra, and G{\"o}ran
  Andersson.
\newblock An uncertainty management framework for integrated gas-electric
  energy systems.
\newblock {\em Proceedings of the IEEE}, 108(9):1518--1540, 2020.

\bibitem{kazi2024intertemporal}
Saif~R. Kazi, Kaarthik Sundar, Sidhant Misra, Svetlana Tokareva, and Anatoly
  Zlotnik.
\newblock Intertemporal uncertainty management in gas-electric energy systems
  using stochastic finite volumes.
\newblock In {\em 2024 Power Systems Computation Conference (PSCC)}. IEEE,
  2024.

\bibitem{geng2021network}
Jiang-Bo Geng, Fu-Rui Chen, Qiang Ji, and Bing-Yue Liu.
\newblock Network connectedness between natural gas markets, uncertainty and
  stock markets.
\newblock {\em Energy Economics}, 95:105001, 2021.

\bibitem{wang2022volatility}
Jiqian Wang, Feng Ma, Elie Bouri, and Juandan Zhong.
\newblock Volatility of clean energy and natural gas, uncertainty indices, and
  global economic conditions.
\newblock {\em Energy Economics}, 108:105904, 2022.

\bibitem{schweppe2013spot}
Fred~C. Schweppe, Michael~C. Caramanis, Richard~D. Tabors, and Roger~E. Bohn.
\newblock {\em Spot pricing of electricity}.
\newblock Springer Science \& Business Media, 2013.

\bibitem{pepper2012implementation}
W.~Pepper, B.~J. Ring, E.~G. Read, and S.~R. Starkey.
\newblock Implementation of a scheduling and pricing model for natural gas.
\newblock {\em Handbook of Networks in Power Systems II}, pages 3--35, 2012.

\bibitem{rudkevich17hicss}
Aleksandr~M. Rudkevich and Anatoly Zlotnik.
\newblock Locational marginal pricing of natural gas subject to engineering
  constraints.
\newblock In {\em Proc. of the 50th Hawaii International Conference on System
  Sciences}, pages 3092--3101, 2017.

\bibitem{zlotnik2019optimal}
Anatoly Zlotnik, Kaarthik Sundar, Aleksandr~M Rudkevich, Aleksandr Beylin, and
  Xindi Li.
\newblock Optimal control for scheduling and pricing intra-day natural gas
  transport on pipeline networks.
\newblock In {\em 2019 IEEE 58th Conference on Decision and Control (CDC)},
  pages 4887--4884. IEEE, 2019.

\bibitem{zlotnik2017economic}
Anatoly Zlotnik, Aleksandr~M. Rudkevich, Evgeniy Goldis, Pablo~A. Ruiz, Michael
  Caramanis, Richard Carter, Scott Backhaus, Richard Tabors, Richard Hornby,
  and Daniel Baldwin.
\newblock Economic optimization of intra-day gas pipeline flow schedules using
  transient flow models.
\newblock In {\em PSIG Annual Meeting}, pages PSIG--1715. PSIG, 2017.

\bibitem{tokareva2024stochastic}
Svetlana Tokareva, Anatoly Zlotnik, and Vitaliy Gyrya.
\newblock Stochastic finite volume method for uncertainty quantification of
  transient flow in gas pipeline networks.
\newblock {\em Applied Mathematical Modelling}, 125:66--84, 2024.

\bibitem{huang2006fitted}
C.-S. Huang, C.-H. Hung, and Song Wang.
\newblock A fitted finite volume method for the valuation of options on assets
  with stochastic volatilities.
\newblock {\em Computing}, 77:297--320, 2006.

\bibitem{herty2010new}
Michael Herty, Jan Mohring, and Veronika Sachers.
\newblock A new model for gas flow in pipe networks.
\newblock {\em Mathematical Methods in the Applied Sciences}, 33(7):845--855,
  2010.

\bibitem{kazi2024modeling}
Saif~R. Kazi, Kaarthik Sundar, Shriram Srinivasan, and Anatoly Zlotnik.
\newblock Modeling and optimization of steady flow of natural gas and hydrogen
  mixtures in pipeline networks.
\newblock {\em International Journal of Hydrogen Energy}, 54:14--24, 2024.

\bibitem{menon05}
E.~S. Menon.
\newblock {\em Gas pipeline hydraulics}.
\newblock CRC Press, 2005.

\bibitem{vuffray2015monotonicity}
Marc Vuffray, Sidhant Misra, and Michael Chertkov.
\newblock Monotonicity of dissipative flow networks renders robust maximum
  profit problem tractable: General analysis and application to natural gas
  flows.
\newblock In {\em 2015 54th IEEE Conference on Decision and Control (CDC)},
  pages 4571--4578. IEEE, 2015.

\bibitem{wu2017adaptive}
Fei Wu, Harsha Nagarajan, Anatoly Zlotnik, Ramteen Sioshansi, and Aleksandr~M
  Rudkevich.
\newblock Adaptive convex relaxations for gas pipeline network optimization.
\newblock In {\em 2017 American Control Conference (ACC)}, pages 4710--4716.
  IEEE, 2017.

\bibitem{sodwatana2023optimization}
Mo~Sodwatana, Saif~R. Kazi, Kaarthik Sundar, and Anatoly Zlotnik.
\newblock Optimization of hydrogen blending in natural gas networks for carbon
  emissions reduction.
\newblock In {\em 2023 American Control Conference (ACC)}, pages 1229--1236.
  IEEE, 2023.

\bibitem{misra2020monotonicity}
Sidhant Misra, Marc Vuffray, and Anatoly Zlotnik.
\newblock Monotonicity properties of physical network flows and application to
  robust optimal allocation.
\newblock {\em Proceedings of the IEEE}, 108(9):1558--1579, 2020.

\bibitem{gyrya2019explicit}
Vitaliy Gyrya and Anatoly Zlotnik.
\newblock An explicit staggered-grid method for numerical simulation of
  large-scale natural gas pipeline networks.
\newblock {\em Applied Mathematical Modelling}, 65:34--51, 2019.

\bibitem{dunning2017jump}
Iain Dunning, Joey Huchette, and Miles Lubin.
\newblock Jump: A modeling language for mathematical optimization.
\newblock {\em SIAM Review}, 59(2):295--320, 2017.

\bibitem{waechter2006ipopt}
Andreas W{\"a}chter and {Lorenz T.} Biegler.
\newblock On the implementation of an interior-point filter line-search
  algorithm for large-scale nonlinear programming.
\newblock {\em Mathematical Programming}, 106(1):25--57, 2006.

\bibitem{srinivasan2022numerical}
Shriram Srinivasan, Kaarthik Sundar, Vitaliy Gyrya, and Anatoly Zlotnik.
\newblock Numerical solution of the steady-state network flow equations for a
  non-ideal gas.
\newblock {\em IEEE Transactions on Control of Network Systems}, 2022.

\end{thebibliography}

\end{document}